\documentclass{amsart}
\usepackage{amssymb}
\usepackage[all]{xy}
\usepackage{amsmath}
\usepackage{latexsym}
\CompileMatrices

\newtheorem{theorem}{Theorem}[section]
\newtheorem{lemma}[theorem]{Lemma}
\newtheorem{proposition}[theorem]{Proposition}
\newtheorem{corollary}[theorem]{Corollary}
\theoremstyle{definition}

\newtheorem{remark}[theorem]{Remark}

\numberwithin{equation}{theorem}

\def\KK{{\mathbb K}}
\def\NN{{\mathbb N}}

\def\codim{{\rm codim}}
\def\Spec{{\rm Spec\,}}
\newcounter{itemnumber}

\begin{document}

\sloppy

\title[Invariant ideals and Matsushima's criterion]
{Invariant ideals and Matsushima's criterion}
\thanks{%
Supported by DFG Schwerpunkt 1094, CRDF grant RM1-2543-MO-03, RF President grant MK-1279.2004.1, and RFBR 05-01-00988.} 
\author[I. Arzhantsev]{Ivan V. Arzhantsev}
\address{Department of Higher Algebra, Faculty of Mechanics and Mathematics, Moscow State Lomonosov University,
Vorobievy Gory, GSP-2, Moscow, 119992, Russia}
\email{arjantse@mccme.ru}
\begin{abstract}
Let $G$ be a reductive algebraic group and $H$ a closed subgroup of $G$. Explicit constructions
of $G$-invariant ideals in the algebra $\KK[G/H]$ are given. This allows to obtain an elementary proof of Matsushima's
criterion: a homogeneous space $G/H$ is an affine variety if and only if $H$ is reductive. 
\end{abstract}

\maketitle


\section{Algebraic homogeneous spaces}

Let $G$ be an affine algebraic group over an algebraically closed field $\KK$. A $G$-module $V$ is said to be rational if
any vector in $V$ is contained in a finite-dimensional rational $G$-submodule. 
Below all modules are supposed to be rational.
By $V^G$ denote the subspace of $G$-fixed vectors in $V$.

The group $G\times G$ acts on $G$ by translations, $(g_1,g_2)g:=g_1gg_2^{-1}$. This action induces the action on the algebra of regular functions
on $G$:
$$
( G\times G) : \KK[G], \ \ ((g_1,g_2)f)(g):=f(g_1^{-1}gg_2).
$$

For any closed subgroup $H$ of $G$, $H_l$, $H_r$ denote the groups of all left and right translations of $\KK[G]$ by 
elements of $H$. 
Under these actions, the algebra $\KK[G]$ becomes a rational $H_l$- (and $H_r$-) module.

\smallskip

By Chevalley's Theorem, the set $G/H$ of left $H$-cosets in $G$ admits a structure of a quasi-projective algebraic variety such that
the projection $p:G\to G/H$ is a surjective $G$-equivariant morphism. Moreover, a structure of an algebraic variety
on $G/H$ satisfying these conditions is unique. It is easy to check that the morphism $p$ is open and the algebra of regular functions on $G/H$ may be
identified with the subalgebra $\KK[G]^{H_r}$ in $\KK[G]$. We refer to~\cite[Ch.~IV]{Hu2} for details.


\section{Matsushima's criterion}

Let $G$ be a reductive algebraic group and $H$ a closed subgroup of $G$. It is known that the homogeneous space $G/H$ is affine
if and only if $H$ if reductive. The first proof was given over the field of complex numbers and used
some results from algebraic topology, see~\cite{Ma} and~\cite[Th.~4]{On}.
An algebraic proof in characteristic zero was obtained in~\cite{BB}. A characteristic-free proof that
uses the Mumford conjecture proved by W.J.~Haboush is given in~\cite{Ri}. Another proof based on the Morozov-Jacobson Theorem may be found
in~\cite{Lu}.   

Below we give an elementary proof of Matsushima's criterion in terms of representation theory. The ground field $\KK$ is assumed to be
algebraically closed and of characteristic zero.


\begin{theorem}\label{t1}

Let $G$ be a reductive algebraic group and $H$ its closed subgroup. Then the homogeneous space $G/H$ is affine if
and only if $H$ is reductive.

\end{theorem}


\begin{proof}

We begin with the ``easy half''.


\begin{proposition}\label{pr1}

Let $G$ be an affine algebraic group and $H$ its reductive subgroup. Then $G/H$ is affine. 

\end{proposition}


\begin{proof}

If a reductive group $H$ acts on an affine variety $X$, then the algebra of invariants $\KK[X]^H$ is finitely
generated, the quotient morphism $\pi: X\to\Spec\KK[X]^H$ is surjective and any fiber of $\pi$ contains
a unique closed $H$-orbit \cite[Sec~.4.4]{PV}. In the case $X=G$ this shows that $G/H$ is isomorphic to $\Spec\KK[G]^{H_r}$.

\end{proof}


Now assume that $G$ is reductive and consider a decomposition
$$
\KK[G]=\KK\oplus\KK[G]_G,
$$

\noindent where the first component corresponds to constant functions on $G$, and the second one is the sum of all
simple non-trivial $G_l$- (or $G_r$-) submodules in $\KK[G]$. Let ${\rm pr}:\KK[G]\to \KK$ be the projection
on the first component. Clearly, ${\rm pr}$ is a $(G_l\times G_r)$-invariant linear map.  

Let $H$ be a closed subgroup of $G$. Consider 

$$
I(G,H)=\{ f\in\KK[G]^{H_r} \ | \ {\rm pr}(fl)=0 \ \forall l\in\KK[G]^{H_r} \}.
$$ 

This is a $G_l$-invariant ideal in $\KK[G]^{H_r}$ with $1\notin I(G,H)$. Assume that $G/H$ is affine. Then 
$G/H\cong\Spec\KK[G]^{H_r}$ and $\KK[G]^{H_r}$ does not contain proper $G_l$-invariant ideals. Thus $I(G,H)=0$. Our aim
is to deduce from this that any $H$-module is completely reducible.


\begin{lemma}\label{l2}

If $W$ is an $H_r$-submodule in $\KK[G]$ and $f\in W$ is a non-zero $H_r$-fixed vector,
then $W=\langle f\rangle\oplus W'$, where $W'$ is an $H_r$-submodule. 

\end{lemma}


\begin{proof}

Since $I(G,H)=0$, there exists $l\in\KK[G]^{H_r}$ such that ${\rm pr}(fl)\ne 0$. 
The submodule $W'$ is defined as $W'=\{ w\in W \ | \ {\rm pr}(wl)=0\}$. 

\end{proof}


\begin{lemma}\label{l1}

If $f\in\KK[G]$ is an $H_r$-semi-invariant of the weight $\xi$, then there exists an $H_r$-semi-invariant in $\KK[G]$ of
the weight $-\xi$.

\end{lemma}


\begin{proof}

Let $Z$ be the zero set of $f$ in $G$. Since $Z$ is $H_r$-invariant, one has $Z=p^{-1}(p(Z))$.
This implies that $p(Z)$ is a proper closed subset of $G/H$. There exists a non-zero $\alpha\in\KK[G/H]$ with $\alpha|_{p(Z)}=0$.
Then $p^*\alpha\in\KK[G]^{H_r}$ and $p^*\alpha|_Z=0$. By Hilbert's Nullstellensatz, there are 
$n\in\NN$, $s\in\KK[G]$ such that $(p^*\alpha)^n=fs$. This shows that $s$ is an $H_r$-semi-invariant of the weight $-\xi$. 

\end{proof}


\begin{lemma}\label{l4}

(1) \ Any cyclic $G$-module $V$ may be embedded (as a $G_r$-submodule) into $\KK[G]$.

(2) \ Any $n$-dimensional $H$-module $W$ may be embedded (as an $H_r$-submodule) into $(\KK[H])^n$.

(3) \ Any finite-dimensional $H$-module may be embedded (as an $H$-submodule) into a finite-dimensional $G$-module.

\end{lemma}


\begin{proof}

(1) Suppose that $V=\langle Gv\rangle$.
The map $\phi: G\to V$, $\phi(g)=g^{-1}v$, induces the embedding of the dual module $\phi^*: V^*\to\KK[G]$. Consider the $G_r$-submodule  
$U=\{f\in\KK[G] \ | \ {\rm pr}(fl)=0 \ \forall l\in\phi^*(V^*)\}$. By the complete reducibility, $\KK[G]=U\oplus U'$ for some
$G_r$-submodule $U'$. Obviously, $I(G,G)=0$ and $U'$ is $G_r$-isomorphic to $V$. 

(2) Let $\lambda_1,\dots,\lambda_n$ be a basis of $W^*$. The embedding may be given as
$$
 w\to (f_1^w,\dots,f_n^w), \ f_i^w(h):=\lambda_i(hw).
$$

(3) Note that the restriction homomorphism $\KK[G]\to\KK[H]$ is surjective. By (2), any finite-dimensional $H$-module
$W$ has the form $W_1/W_2$, where $W_1$ is a finite-dimensional  $H$-submodule in a $G$-module $V$ and $W_2$ is an $H$-submodule of $W_1$.
Consider $W_1\bigwedge(\bigwedge^m W_2)$ as an $H$-submodule in $\bigwedge^{m+1} W_1$, where $m=\dim W_2$. Note that
$W\cong (W_1\bigwedge(\bigwedge^m W_2))\otimes (\bigwedge^m W_2)^*$. By (1), the cyclic $G$-submodule of $\bigwedge^m V$ generated
by $\bigwedge^m W_2$ may be embedded into $\KK[G]$. By Lemma~\ref{l1}, $(\bigwedge^m W_2)^*$ also may be embedded into a 
$G$-module.
 
\end{proof}


\begin{lemma}\label{l3}

For any $H$-module $W$ and any non-zero $w\in W^H$ there is an $H$-submodule $W'$
such that $W=\langle w\rangle\oplus W'$.

\end{lemma}


\begin{proof}

Embed $W$ into a $G$-module $V$. Let $V_1=\langle Gw\rangle$. Then $V=V_1\oplus V_2$ for some $G$-submodule $V_2$.
Embed $V_1$ into $\KK[G]$ as a $G_r$-submodule. By Lemma~\ref{l2}, $V_1=\langle w\rangle\oplus W_1$ for some
$H$-submodule $W_1$. Finally, $W'=W\cap (W_1\oplus V_2)$.

\end{proof}


\begin{lemma}\label{wl}

Any $H$-module is completely reducible.

\end{lemma}


\begin{proof}

Assume that $W_1$ is a simple submodule in an $H$-module $W$.
Consider two submodules in the $H$-module ${\rm End}(W,W_1)$:

$$
  L_2=\{ p\in {\rm End}(W,W_1) \ | \ p|_{W_1}=0\} \ \subset \ L_1=\{p\in{\rm End}(W,W_1) \ | \ p|_{W_1} \ {\rm is \ scalar} \}.
$$

Clearly, $L_2$ is a hyperplane in $L_1$. Consider an $H$-eigenvector $l\in (L_1)^*$ corresponding to $L_2$. Taking the tensor product
with a one-dimensional $H$-module, one may assume that $l$ is $H$-fixed. By Lemma~\ref{l3}, $(L_1)^*=\langle l\rangle\oplus M$, implying
$L_1=L_2\oplus\langle P\rangle$, where $M$ and $\langle P\rangle $ are $H$-submodules. 
Then ${\rm Ker}\,P$ is a complementary submodule to $W_1$.    

\end{proof}


Theorem~\ref{t1} is proved.

\end{proof}


\begin{remark}

In~\cite{Vi}, for any action $G:X$ of a reductive group $G$ on an affine variety $X$ with the decomposition
$\KK[X]=\KK[X]^G\oplus\KK[X]_G$ and the projection ${\rm pr}:\KK[X]\to\KK[X]^G$, the $\KK[X]^G$-bilinear scalar
product $(f,g)={\rm pr}(fg)$ on $\KK[X]$ was introduced and the kernel of this product was considered.
Our ideal $I(G,H)$ is such kernel in the case $X=\Spec\KK[G]^{H_r}$ provided $\KK[G]^{H_r}$ is finitely
generated. 

\end{remark}


\begin{remark}

For convenience of the reader we include all details in the proof of Theorem~\ref{t1}. Lemma~\ref{l1} and Lemma~\ref{l4} are
taken from~\cite{BBHM}. They show that for a quasi-affine $G/H$ any $H$-module may be realized as an $H$-submodule of a $G$-module.
The converse is also true~\cite{BBHM},~\cite{Gr}. 
Proposition~\ref{pr1} is a standart fact. The proof of Lemma~\ref{wl} is a part of the proof of the Weyl Theorem
on complete reducibility~\cite{Hu}, see also~\cite[Prop.~2.2.4]{Sp}.

\end{remark}


\section{Some additional remarks}

The following lemma may be found in~\cite{BB}.

\begin{lemma}\label{lbn}

Let $G$ be an affine algebraic group and $H$ its reductive subgroup. Then $\KK[G]^{H_r}$ does not contain
proper $G_l$-invariant ideals.

\end{lemma}


\begin{proof}

Consider a decomposition
$$
 \KK[G]=\KK[G]^{H_r}\oplus\KK[G]_{H_r},
$$

\noindent where $\KK[G]_{H_r}$ is the sum of all non-trivial simple $H_r$-submodules in $\KK[G]$. Clearly, 
$\KK[G]^{H_r}\KK[G]_{H_r}\subseteq\KK[G]_{H_r}$. Hence any proper $G_l$-invariant ideal in $\KK[G]^{H_r}$ 
generates a proper $G_l$-invariant ideal in $\KK[G]$, a contradiction.

\end{proof}


By Hilbert's Theorem on invariants, the algebra $\KK[G]^{H_r}$ is finitely generated.
It is easy to see that functions from $\KK[G]^{H_r}$ separate (closed) right $H$-cosets in $G$. These observations 
and Lemma~\ref{lbn} give another proof of
Proposition~\ref{pr1}. Moreover, it is proved in~\cite[Prop.~1]{BB} that for a quasi-affine $G/H$ the algebra
$\KK[G]^{H_r}$ does not contain proper $G_l$-invariant ideals if and only if $G/H$ is affine. 


\smallskip

Now assume that $G$ is reductive.

\begin{proposition}\cite[Prop.~1]{Vi}\label{pr2}
The ideal $I(G,H)$ is the biggest $G_l$-invariant ideal in $\KK[G]^{H_r}$ different from $\KK[G]^{H_r}$.

\end{proposition}


\begin{proof}

Any proper $G_l$-invariant ideal $I$ of $\KK[G]^{H_r}$ is contained in $\KK[G]^{H_r}\cap\KK[G]_G$. 
Thus ${\rm pr}(il)=0$ for any $l\in\KK[G]^{H_r}$, $i\in I$. This implies $I\subseteq I(G,H)$.

\end{proof}


\begin{remark}

For non-reductive $G$ the biggest invariant ideal in $\KK[G]^{H_r}$ may not exist. For example,
one may take 
$$
G=\left\{
\begin{pmatrix}
1 & * & * \\
0 & * & * \\
0 & * & * 
\end{pmatrix}
\right\}, \ \ 
H=\left\{
\begin{pmatrix}
1 & 0 & * \\
0 & 1 & * \\
0 & 0 & *
\end{pmatrix}
\right\}.
$$

Here $G/H\cong \KK^3\setminus\{x_2=x_3=0\}$, $\KK[G]^{H_r}\cong\KK[x_1,x_2,x_3]$, and the maximal ideals
$(x_1-a,x_2,x_3)$ are $G_l$-invariant for any $a\in\KK$. 

\end{remark}


\section{The boundary ideal}

In this section we assume that $H$ is an observable subgroup of $G$, i.e., $G/H$ is quasi-affine. 

If the algebra $\KK[G]^{H_r}$ is finitely generated, then the affine $G$-variety $X=\Spec\KK[G]^{H_r}$ has an
open $G$-orbit isomorphic to $G/H$ and may be considered as the canonical embedding $G/H\hookrightarrow X$. Moreover,
this embedding is uniquely characterized by two properties: $X$ is normal and $\codim_X (X\setminus G/H)\ge 2$, see~\cite{Gr}.
There are two remarkable $G_l$-invariant ideals in $\KK[G]^{H_r}$, namely
$$
 I^b(G,H)=I(X\setminus(G/H))=\{f\in\KK[G]^{H_r} \ | \ f|_{X\setminus(G/H)}=0\},
$$
and, if $G$ is reductive, the ideal $I^m(G,H)$ of the unique closed $G$-orbit in $X$. If $G/H$ is affine, then $I^b(G,H)=\KK[G]^{H_r}$,
$I^m(G,H)=0$. In other cases $I^b(G,H)$ is the smallest proper radical $G_l$-invariant ideal, and $I^m(G.H)$ is the biggest
proper $G_l$-invariant ideal of $\KK[G]^{H_r}$. 
By Proposition~\ref{pr2}, $I^m(G,H)=I(G,H)$. Moreover, 
$\KK[G]^{H_r}/I^m(G,H)\cong\KK[G]^{S_r}$, where $S$ is a minimal reductive subgroup of $G$ containing $H$. (Such a subgroup may be not unique, 
but all of them are $G$-conjugate, see~\cite[Sec.~7]{Ar}.) 
It follows from the Slice Theorem~\cite{Lu} and \cite[Prop.~4]{Ar} that $I^b(G,H)=I^m(G,H)$ if and only if $H$ is a
quasi-parabolic subgroup of a reductive subgroup of $G$.

\smallskip

Now assume that $\KK[G]^{H_r}$ is not finitely generated. If $G$ is reductive, then $I(G,H)$ may be consider as an analog of $I^m(G,H)$
in this situation (Proposition~\ref{pr2}). We claim that $I^b(G,H)$ also has an analog, even for non-reductive $G$. 


\begin{proposition}\label{pr3}

Let $\hat X$ be a quasi-affine variety, $\hat X\hookrightarrow X$ be an (open) embedding into an affine variety $X$,
$I(X\setminus \hat X)\lhd\KK[X]$, and $\mathcal{I}=\mathcal{I}(\hat X)$ be the radical of the ideal of $\KK[\hat X]$ generated by
$I(X\setminus\hat X)$. Then

(1) \  the ideal $\mathcal{I}\lhd\KK[\hat X]$ does not depend on $X$;

(2) \ $I(X\setminus\hat X)$ is the smallest radical ideal of $\KK[X]$ generating an ideal in $\KK[\hat X]$ with the radical $\mathcal{I}$.
 
\end{proposition}
 

\begin{proof}

(1) \ Consider two affine embeddings: $\phi_i:\hat X\hookrightarrow X_i$, $i=1,2$. Let $X_{12}$ be the closure of $(\phi_1\times\phi_2)(\hat X)$ in
$X_1\times X_2$ with the projections $r_i: X_{12}\to X_i$. Let us identify the images of $\hat X$ in $X_1$, $X_2$, and $X_{12}$. We claim
that $r_i(X_{12}\setminus\hat X)\subseteq X_i\setminus\hat X$. Indeed, the diagonal image of $\hat X$ is closed in $\hat X\times X_j$, $j\ne i$, as
the graph of a morphism.


It follows from what was proved above that the ideal of $\KK[X_{12}]$ generated by $r_i^*(I(X_i\setminus\hat X))$ has the radical
$I(X_{12}\setminus\hat X)$. This shows that the radical of the ideal generated by $I(X_i\setminus\hat X)$ in $\KK[\hat X]$ does not
depend on $i$.  

(2) Assume that there is a radical ideal $I_1\lhd\KK[X]$ not containing $I=I(X\setminus\hat X)$ and generating 
an ideal in $\KK[\hat X]$ with the radical $\mathcal{I}$.
There is $x_0\in\hat X$ such that $h(x_0)=0$ for any $h\in I_1$. Take $f\in I$ such that $f(x_0)\ne 0$.
One has
$f^k=\alpha_1h_1+\dots+\alpha_kh_k$ for some $\alpha_i\in\KK[\hat X]$, $h_i\in I_1$, $k\in\NN$, and this implies $f(x_0)=0$, a contradiction. 

\end{proof}


So $\mathcal{I}(G/H)$ is a radical $G_l$-invariant ideal of $\KK[G]^{H_r}$, and $\mathcal{I}(G/H)=I^b(G,H)$ provided
$\KK[G]^{H_r}$ is finitely generated.


\begin{proposition}

$\mathcal{I}(G/H)$ is the smallest non-zero radical $G_l$-invariant ideal of $\KK[G]^{H_r}$.

\end{proposition}


\begin{proof}

Let $f\in\KK[G]^{H_r}$ and $I(f)$ be the ideal of $\KK[G]^{H_r}$ generated by the orbit $G_lf$. It is sufficient
to prove that $\mathcal{I}(G/H)\subseteq {\rm rad}\,I(f)$. Take any $G$-equivariant affine embedding
$G/H\hookrightarrow X$ with $f\in\KK[X]$. For the ideal $I'(f)$ generated by $G_lf$ in $\KK[X]$
one has $I(X\setminus(G/H))\subseteq {\rm rad}\,I'(f)$, hence $\mathcal{I}(G/H)\subseteq{\rm rad}\,I(f)$.

\end{proof}


\begin{corollary}

Let $G$ be an affine algebraic group and $H$ its observable subgroup. Then $G/H$ is affine if and only if
$\mathcal{I}(G/H)=\KK[G]^{H_r}$.

\end{corollary}


It should be interesting to give a description of the ideal $\mathcal{I}(G/H)$ similar to the 
definition of $I(G,H)$,
and to find a geometric meaning of the $G_l$-algebras $\KK[G]^{H_r}/I(G,H)$ and $\KK[G]^{H_r}/\mathcal{I}(G/H)$
for non-finitely generated $\KK[G]^{H_r}$.

\medskip

{\it Acknowledgements.} The author is grateful to J.~Hausen for useful discussions. In particular, Proposition~\ref{pr3}
appears during such a discussion. Thanks are also due to D.A.~Timashev for valuable remarks.

This paper was written during the staying at Eberhard Karls Universit\"at T\"ubingen (Germany). The author wishes to
thank this institution and especially J\"urgen Hausen for invitation and hospitality. 



\begin{thebibliography}{}

\bibitem{Ar} I.V.~Arzhantsev, {\it Algebras with finitely generated invariant subalgebras},
Ann. Inst. Fourier {\bf 53:2} (2003), 379-398.


\bibitem{BB} A.~Bialynicki-Birula, {\it On homogeneous affine spaces of linear algebraic groups},
Amer. J. Math. {\bf 85} (1963), 577--582. 

\bibitem{BBHM} A.~Bialynicki-Birula, G.~Hochschild and G.D.~Mostow {\it Extensions of representations of algebraic linear groups},
Amer. J. Math. {\bf 85} (1963), 131--144. 

\bibitem{Gr} F.D.~Grosshans. Algebraic Homogeneous Spaces and Invariant Theory.
LNM {\bf 1673}, Springer-Verlag, Berlin, 1997.

\bibitem{Hu} J.E.~Humphreys. Introduction to Lie Algebras and Representation Theory.
Graduate Texts in Math. {\bf 9}, Springer-Verlag, New York Berlin, 1972.

\bibitem{Hu2} J.E.~Humphreys. Linear Algebraic Groups.
Graduate Texts in Math. {\bf 21}, Springer-Verlag, New York Heidelberg Berlin, 1975.

\bibitem{Lu} D.~Luna, {\it Slices \'etales}, Bull. Soc. Math. Fr., 
memoire {\bf 33} (1973), 81--105.

\bibitem{Ma} Y.~Matsushima, {\it Espaces homog\'enes de Stein des groupes de Lie complexes},
Nagoya Math. J. {\bf 16} (1960), 205--218.

\bibitem{On} A.L.~Onishchik, {\it Complex hulls of compact homogeneous spaces},
Dokl. Akad. Nauk SSSR {\bf 130:4} (1960), 726--729 (Russian);
English Transl.: Soviet Math. Dokl {\bf 1} (1960), 88--91.

\bibitem{PV} V.L.~Popov and E.B.~Vinberg. Invariant Theory. Itogi Nauki i Tekhniki, Sovrem.
Probl. Mat. Fund. Naprav. vol.~{\bf 55}, VINITI, Moscow, 1989. 
English Transl.: Algebraic Geometry IV, Encyclopedia of Math. Sciences, vol.~{\bf 55},
Springer-Verlag, Berlin, 1994.

\bibitem{Ri} R.W.~Richardson, {\it Affine coset spaces of reductive algebraic groups},
Bull. London Math. Soc. {\bf 9:1} (1977), 38--41.

\bibitem{Sp} T.A.~Springer. Invariant Theory. LNM {\bf 585}, Springer-Verlag,
Berlin Heidelberg New York, 1977.

\bibitem{Vi} E.B.~Vinberg, {\it On stability of actions of reductive algebraic groups},
in ``Lie Algebras, Rings and Related Topics'', Fong Yuen, A.A.~Mikhalev, E.~Zelmanov Eds.,
Springer-Verlag, Hong Kong Ltd. (2000), 188--202.

\end{thebibliography}
\end{document}